%% file: Art-E2VEM-Comparison.tex
\documentclass[11pt,a4paper]{article}
\usepackage{Art-E2VEM-Comparison}

\input{frontmatter}

\begin{document}
\maketitle

\begin{abstract}
  In this letter we compare the behaviour of standard Virtual Element Methods
  (VEM) and stabilization free Enlarged Enhancement Virtual Element Methods
  (\EEVEM{}) with the focus on some elliptic test problems whose solution and
  diffusivity tensor are characterized by anisotropies. Results show that the
  possibility to avoid an arbitrary stabilizing part, offered by \EEVEM{}
  methods, can reduce the magnitude of the error on general polygonal meshes and
  help convergence.
\end{abstract}


\input{model_problem}
\input{discrete_problem}
\input{numerical_results}
\input{conclusions}

\bibliographystyle{plain}
\bibliography{bibliografia}

\end{document}

%% file: frontmatter.tex
  \title{Comparison of standard and stabilization free Virtual Elements on anisotropic elliptic
    problems}

  \author{Stefano Berrone, Andrea Borio, Francesca
    Marcon\footnote{The three authors
    are members of the Gruppo Nazionale Calcolo Scientifico (GNCS) at Istituto Nazionale di Alta
    Matematica (INdAM). The authors kindly acknowledge partial financial support by INdAM-GNCS
    Projects 2020, by the Italian Ministry of Education, University and Research (MIUR) through the
    ``Dipartimenti di Eccellenza'' Programme (2018–2022) – Department of Mathematical Sciences
    ``G. L. Lagrange'', Politecnico di Torino (CUP:E11G18000350001) and through the PRIN 2017
    project (No. 201744KLJL\_004).}
    } 
    

%% file: model_problem.tex
\section{Introduction}
\label{sec:model_problem}
In recent years, polytopal methods for the solution of PDEs have received a huge
attention from the scientific community. VEM were introduced in
\cite{Beirao2013a,Ahmad2013,Beirao2015b} as a family of methods that deal with
polygonal and polyhedral meshes without building an explicit basis of functions
on each element, but rather defining the local discrete spaces and degrees of
freedom in such a way that suitable polynomial projections of basis functions
are computable. The problem is discretized with bilinear forms that consist of a
polynomial part that mimics the operator and an arbitrary stabilizing bilinear
form. In \cite{ABBDVW}, error analysis focused on anisotropic elliptic problems
shows that the stabilization term adds an isotropic component of the error,
independently of the nature of the problem. In \cite{BBME2VEM}, a modified
version of the method, \EEVEM{}, was proposed, designed to allow the definition
of coercive bilinear forms that consist only of a polynomial approximation of
the problem operator. In this letter, we apply the two methods to solve some
test Laplace problems
with anisotropic solutions and diffusivity
tensors. For each test, we compare the relative energy errors done by each
method.

Let $\Omega \subset \mathbb{R}^2$ be a bounded open set with Lipschitz
boundary. We look for a solution of Laplace problem with homogeneous Dirichlet
boundary conditions, that in variational form reads: find
$u\in\sobh[0]{1}{\Omega}$ such that
\begin{equation}
  \label{eq:model}
  \scal[\Omega]{\K\nabla u}{\nabla v} = \scal[\Omega]{f}{v} \quad
  \forall v \in \sobh[0]{1}{\Omega} \,,
\end{equation}
where $\scal[\Omega]{\cdot}{\cdot}$ denotes the $\lebl{\Omega}$ scalar product
and we assume $f\in\lebl{\Omega}$ and
$\K \in [\lebl[\infty]{\Omega}]^{2\times 2}$ is a symmetric positive definite
matrix.


%% file: discrete_problem.tex
\section{Problem discretization}
\label{sec:discr_prob}
We consider a star-shaped polygonal tessellation $\Mh$ of $\Omega$ satisfying
the standard VEM regularity assumptions (see \cite{Beirao2015b,BBME2VEM}).  Let
$k\in\mathbb{N}$ such that $k\geq 1$ and, $\forall E \in\Mh$, let
$\proj[\nabla]{k}{E}\colon \sobh{1}{E} \to \Poly{k}{E}$ be such that,
$\forall v\in\sobh{1}{E}$,
\begin{equation*}
      \scal[E]{\nabla \proj[\nabla]{k}{E}v}{\nabla p} = \scal[E]{\nabla v}{\nabla
      p} \; \forall p \in\Poly{k}{E} \;\text{and}\;
  \begin{cases}
    \int_{\partial E}\proj[\nabla]{k}{E}v = \int_{\partial E}v &
    \text{if $k=1$}\,,
    \\
    \int_E\proj[\nabla]{k}{E}v = \int_Ev & \text{if $k>1$}\,.
  \end{cases}
\end{equation*}

\subsection{Standard Virtual Element discretization}
\label{sec:standard-vem}

According to \cite{Beirao2015b}, we define the following virtual space on any
$E\in\Mh$:
\begin{multline}
  \label{eq:defstandard-vem}
  \Vh[E] = \left\{ v_h\in\sobh{1}{E}\colon \Delta v_h\in \Poly{k}{E},\,
    \eval{v_h}{e} \in\Poly{k}{e} \;\forall e\subset \partial E,\,
    \eval{v_h}{\partial E}\in\cont{\partial E},\right.
  \\
  \left.  \scal[E]{v_h}{p} = \scal[E]{\proj[\nabla]{k}{E}v}{p} \; \forall
    p\in\Poly{k}{E}/\Poly{k-2}{E} \right\} \,,
\end{multline}
and the relative global space
$\Vh = \{v_h\in\sobh[0]{1}{\Omega}\colon \eval{v_h}{E}\in\Vh[E]\; \forall E
\in\Mh\}$. Then \eqref{eq:model} can be discretized by defining,
$\forall E \in\Mh$, the stabilizing bilinear
$\vemstab[E]{}{}\colon\Vh[E]\times\Vh[E]\to \mathbb{R}$ such that, denoting by
$\chi^E(v_h)$ the vector of degrees of freedom of $v_h$ on $E$ (see
\cite{Beirao2015b}),
\begin{equation*}
  \vemstab[E]{u_h}{v_h} = \chi^E(u_h)\cdot \chi^E(v_h) \quad \forall u_h, v_h\in\Vh[E] \,,
\end{equation*}
and looking for $u_h^{\mathcal{V}}\in \Vh$ that solves
\begin{multline}
  \label{eq:discr:vemstandard}
  \sum_{E\in\Mh}\scal[E]{\K\proj{k-1}{E}\nabla u_h^{\mathcal{V}}}{\proj{k-1}{E}\nabla v_h} +
  \norm[{\lebl[\infty]{E}}]{\K}\vemstab[E]{(I-\proj[\nabla]{k}{E})u_h^{\mathcal{V}}}
  {(I-\proj[\nabla]{k}{E})v_h}
  \\
  = \sum_{E\in\Mh}\scal[E]{f}{\proj{k-1}{E}v_h} \quad \forall v_h\in\Vh\,,
\end{multline}
where $\proj{k-1}{E}$ denotes the $\lebl{E}$-projection on $\Poly{k-1}{E}$ or
$\PolyDouble{k-1}{E}$, depending on the context.

\subsection{Enlarged Enhancement Virtual Element discretization}
\label{sec:e2vem}

In \cite{BBME2VEM}, the space defined in \eqref{eq:defstandard-vem} has been
modified in order to allow the discrete problem to be well-posed without the
need of defining a stabilizing bilinear form. Let $\ell_E\in\mathbb{N}$ be given
$\forall E$, such that, denoting by $N_E$ the number of vertices of $E$,
\begin{equation*}
  (k+\ell_E)(k+\ell_E+1) \geq k N_E + k(k+1) -3\,.
\end{equation*}
We define
\begin{multline}
  \label{eq:defe2vem}
  \Wh[E] = \left\{ v_h\in\sobh{1}{E}\colon \Delta v_h\in \Poly{k+\ell_E}{E},\,
    \eval{v_h}{e} \in\Poly{k}{e} \;\forall e\subset \partial E,\, \right.
  \\
  \left.\eval{v_h}{\partial E}\in\cont{\partial E},\, \scal[E]{v_h}{p} =
    \scal[E]{\proj[\nabla]{k}{E}v}{p} \; \forall
    p\in\Poly{k+\ell_E}{E}/\Poly{k-2}{E} \right\} \,,
\end{multline}
that can be seen to have the same degrees of freedom of $\Vh[E]$. Let
$\Wh = \{v_h \in \sobh[0]{1}{\Omega}\colon \eval{v_h}{E}\in\Wh[E]\;\forall E \in
\Mh\}$. Then, we can discretize \eqref{eq:model} by looking for $u_h^{\mathcal{W}}\in\Wh$
such that, $\forall v_h\in\Wh$,
\begin{equation}
  \label{eq:discr_e2vem}
  \sum_{E\in\Mh}\scal[E]{\K\proj{k+\ell_E-1}{E}\nabla u_h^{\mathcal{W}}}
  {\proj{k+\ell_E-1}{E}\nabla v_h} =
  \sum_{E\in\Mh}\scal[E]{f}{\proj{k-1}{E}v_h} \,.
\end{equation}
The proof of well-posedness of \eqref{eq:discr_e2vem} for $k=1$ can be found in
\cite{BBME2VEM}, while its extension to $k>1$ will be the subject of an upcoming
work.

%% file: numerical_results.tex
\section{Numerical results}
\begin{figure}
  \centering
  \begin{subfigure}{.49\linewidth}
    \centering \includegraphics[width=.5\linewidth]{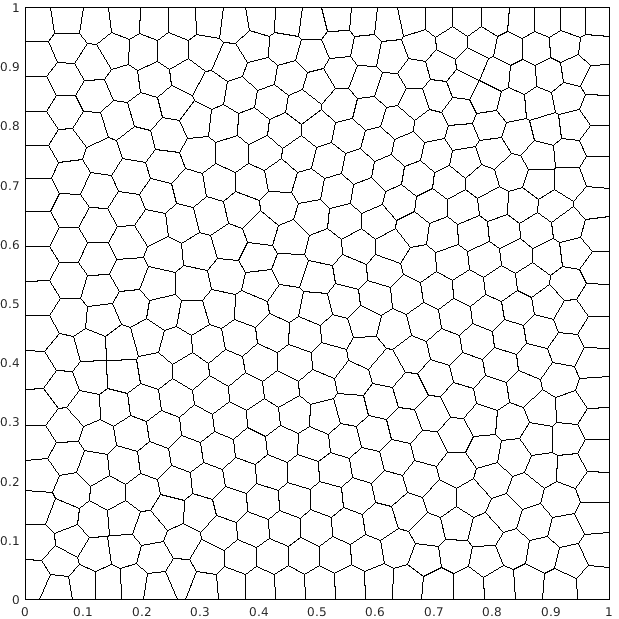}
    \caption{\texttt{Polymesher}}
    \label{fig:meshes:poly}
  \end{subfigure}
  \begin{subfigure}{.49\linewidth}
    \centering \includegraphics[width=.5\linewidth]{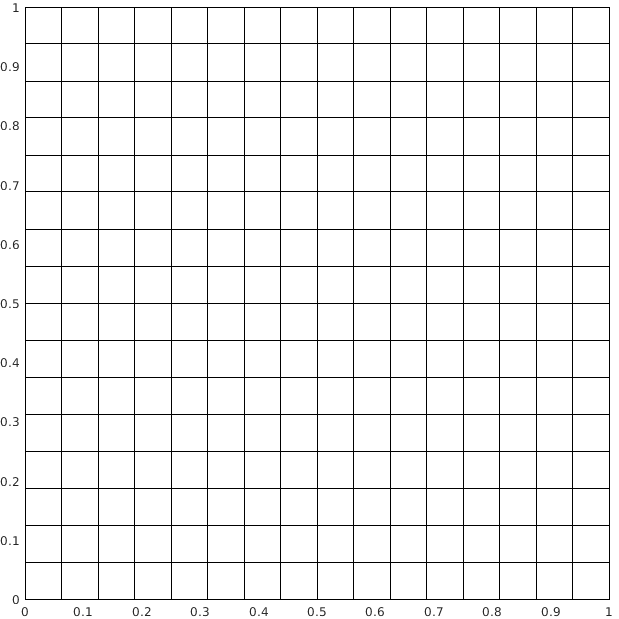}
    \caption{\texttt{Cartesian}}
    \label{fig:meshes:cart}
  \end{subfigure}
  \caption{Meshes used in tests.}
  \label{fig:meshes}
\end{figure}

In all the test cases, we consider problem \eqref{eq:model} on the unit
square. We discretize the domain with the two families of polygonal meshes that
are depicted in Figure \ref{fig:meshes}, the first one being obtained using
Polymesher \cite{Polymesher}, while the second one is a family of standard
cartesian meshes. We compare the two methods described in the previous section
by observing the behaviour of the relative error computed in energy norm as
\begin{equation*}
  e^\star = \frac{\left(\sum_{E\in\Mh}\norm[\lebl{E}]{\sqrt{\K}\nabla\left(u - \proj[\nabla]{k}{E}u_h^\star\right)}^2\right)^{\frac12}}{\norm[\lebl{\Omega}]{\sqrt{\K}\nabla u}} \quad \star = {\mathcal{V}},{\mathcal{W}}\,.
\end{equation*}
In the plots, we show the rate of convergence $\alpha$ computed using the last
two computed errors.

\subsection{Test case 1}
\begin{figure}
  \centering
  \begin{subfigure}{.49\linewidth}
    \centering
    \includegraphics[width=.8\linewidth]{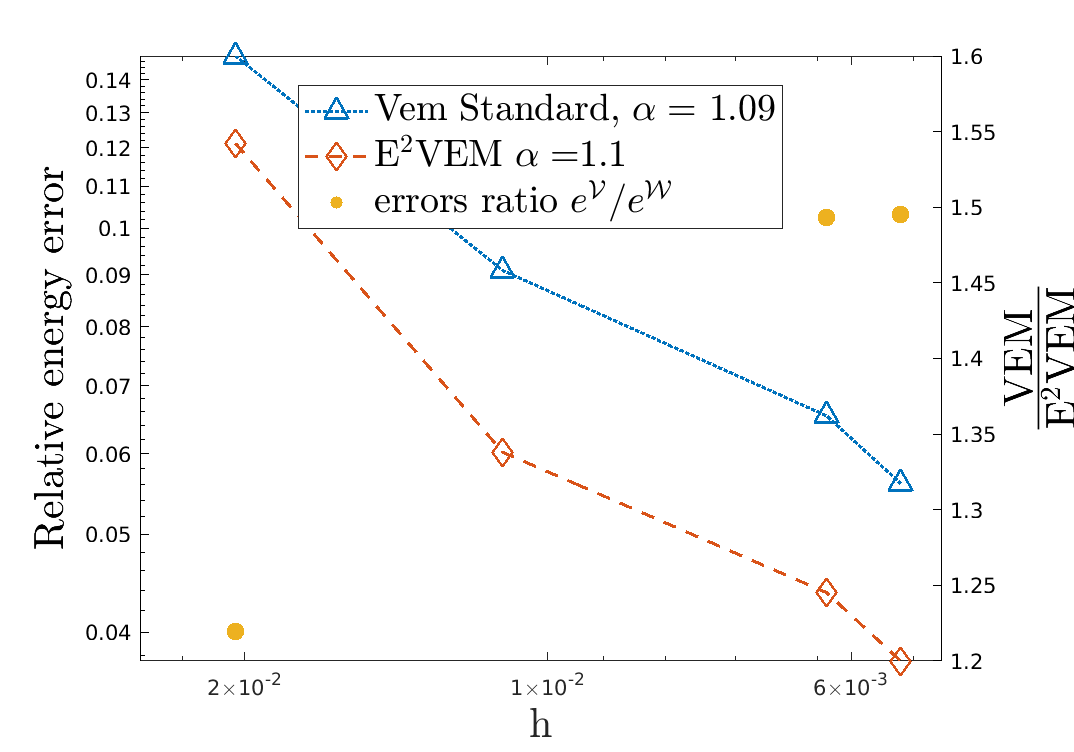}
    \caption{Order 1, \texttt{Polymesher} mesh}
    \label{fig:test2:order1:polymesher}
  \end{subfigure}
  \begin{subfigure}{.49\linewidth}
    \centering
    \includegraphics[width=.8\linewidth]{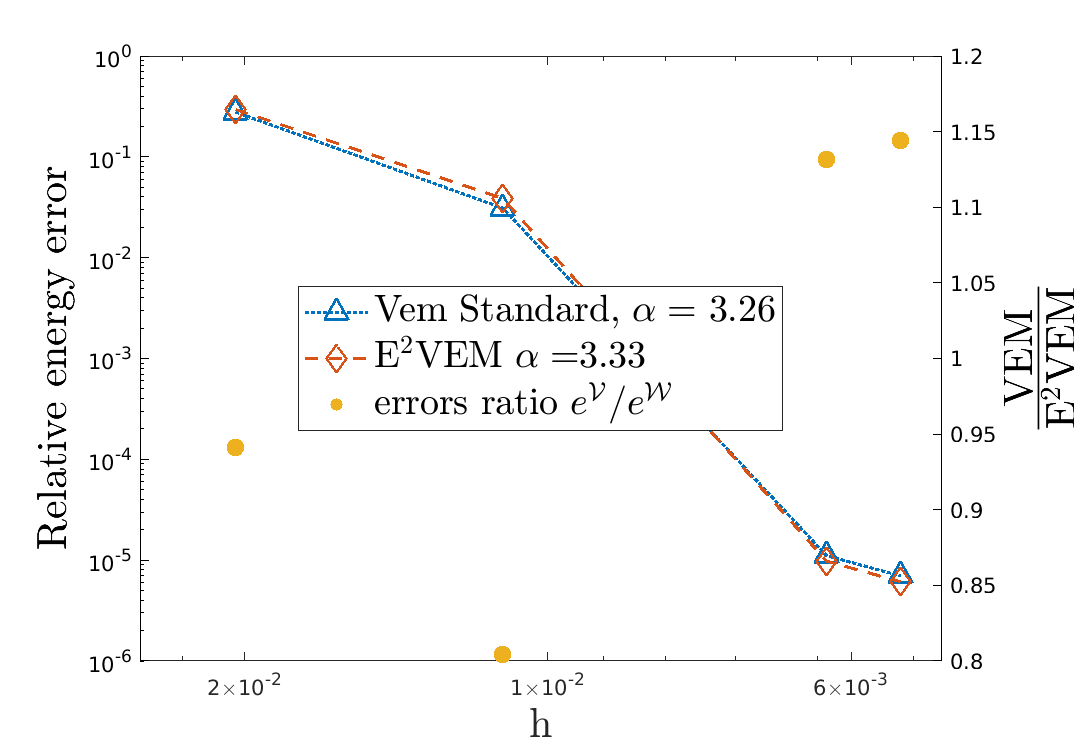}
    \caption{Order 3, \texttt{Polymesher} mesh}
  \end{subfigure}
  \begin{subfigure}{.49\linewidth}
    \centering
    \includegraphics[width=.8\linewidth]{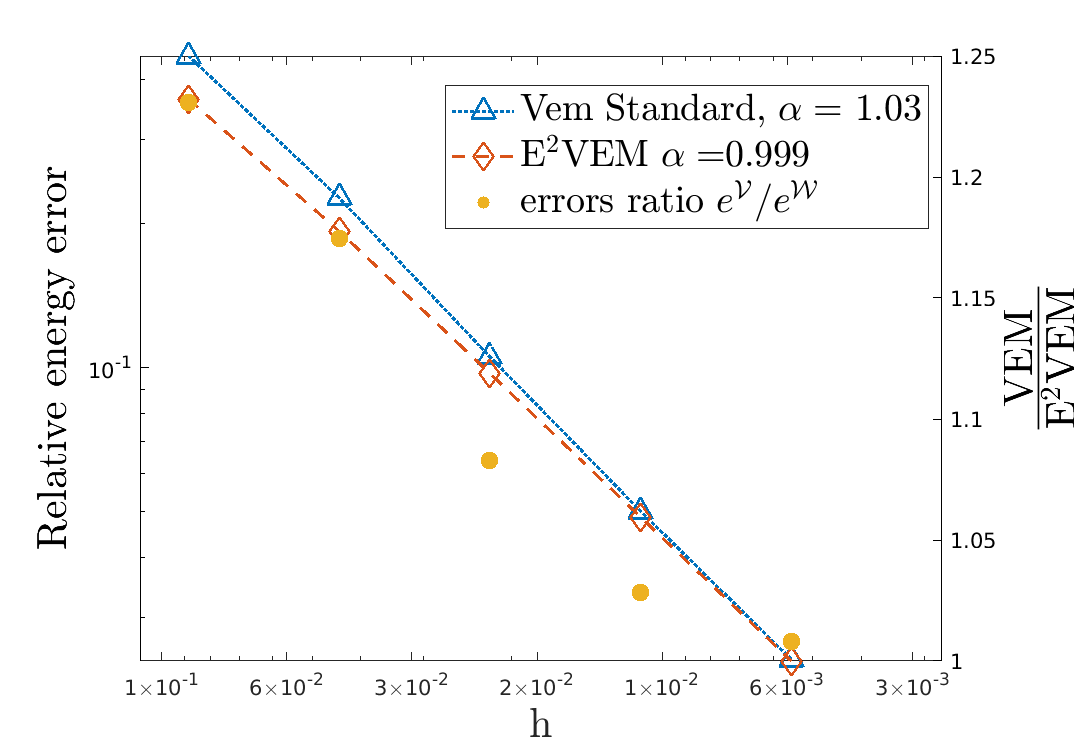}
    \caption{Order 1, \texttt{Cartesian} mesh}
  \end{subfigure}
  \begin{subfigure}{.49\linewidth}
    \centering
    \includegraphics[width=.8\linewidth]{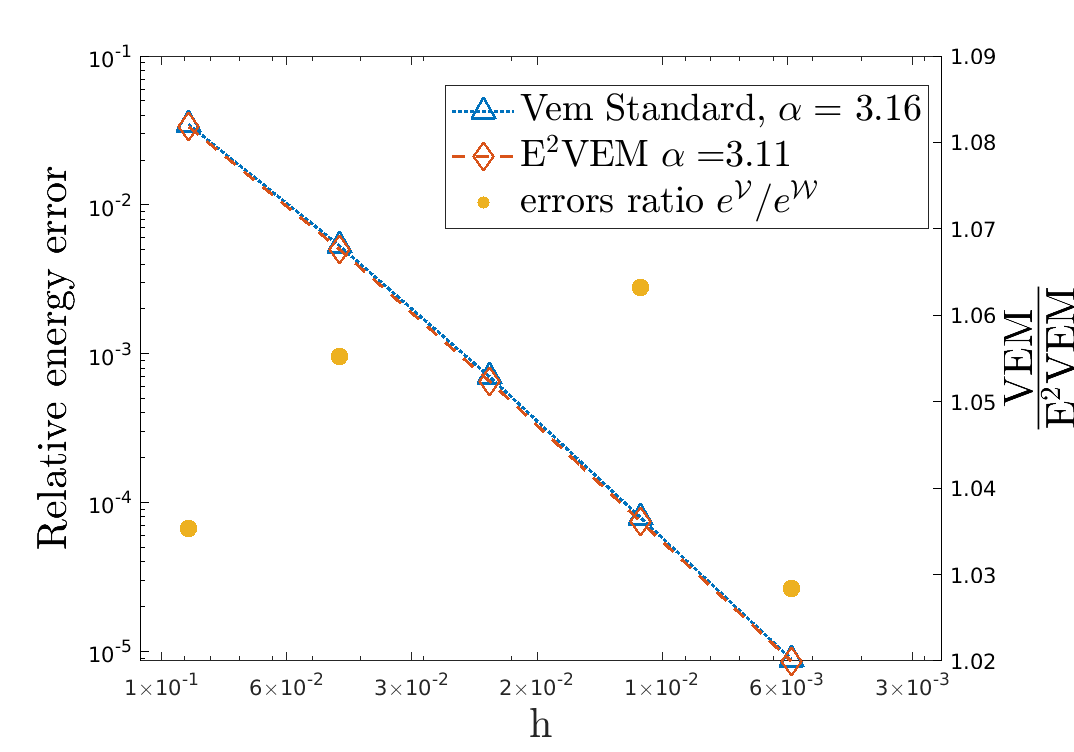}
    \caption{Order 3, \texttt{Cartesian} mesh}
  \end{subfigure}
  \caption{Test case 1. Convergence plots.}
  \label{fig:test2}
\end{figure}
\begin{table}
  \centering \input{tab/test2_projstabratio.tex}
  \caption{Test case 1. Average ratio through refinement between the infinity
    norms of the polynomial part $\mathbb{A}^\Pi$ and the stabilizing part
    $\mathbb{A}^S$ of the stiffness matrix in standard VEM.}
  \label{tab:test2_projstabratio}
\end{table}
In the first test, we define the forcing term $f$ such that the exact solution is
$u(x,y)=10^{-2}xy(1-x)(1-y)(\mathrm{e}^{20x}-1)$, whereas
$\K = 8\cdot10^{-3}(e_1e_1^\intercal) + e_2e_2^\intercal $, where $e_1$ and
$e_2$ are the vectors of the canonical basis of $\mathbb{R}^2$. Figure
\ref{fig:test2} displays the behaviour of the errors obtained with the two methods and the ratio $e^{\mathcal{V}}/e^{\mathcal{W}}$, for orders 1
 and 3, with respect to the maximum diameter of the discretization. The results show that the two methods behave equivalently on
cartesian meshes, whereas \EEVEM{} performs better on the \texttt{Polymesher}
meshes with order $1$, while the two methods tend to have the same behaviour
with higher orders. This is due to the strong anisotropy both of the solution
(that presents a strong boundary layer in the $x$-direction close to the
boundary $x=1$ of the domain) and of the diffusivity tensor $\K$. Indeed, as we
can see from Table \ref{tab:test2_projstabratio}, for $k=1$ the stabilizing part
of the VEM bilinear form is of the same order of magnitude as the polynomial
part, while for $k=3$ we can see that the polynomial part is predominant. This
induces larger errors (see Figure \ref{fig:test2:order1:polymesher}) for the
standard method on general polygonal meshes, such as the ones in the
\texttt{Polymesher} family, since the stabilization is an isotropic
operator. This effect is not felt by the \EEVEM{} method since its bilinear form
consists only on a polynomial part that correctly takes into account the
anisotropy of the tensor $\K$. The difference between the two methods is
mitigated on \texttt{Cartesian} meshes since they are by construction aligned
with the principal directions of the error (see the error analysis done in
\cite{ABBDVW}).

\subsection{Test case 2}
\begin{figure}
  \centering
  \begin{subfigure}{.49\linewidth}
    \centering
    \includegraphics[width=.8\linewidth]{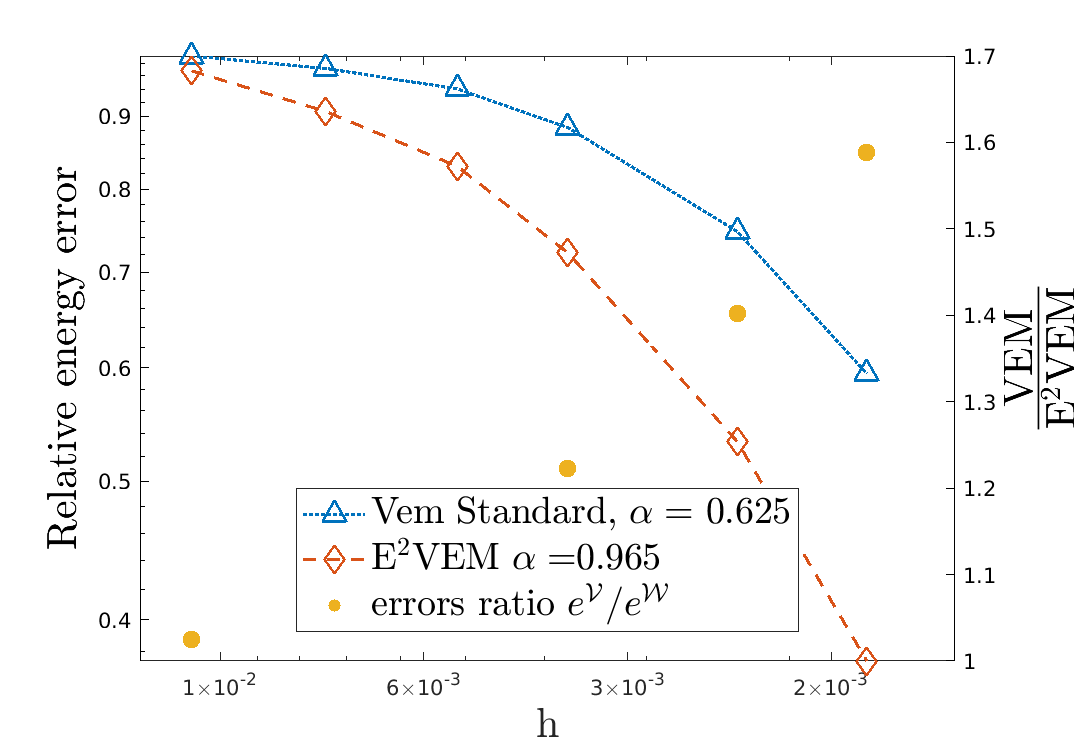}
    \caption{Order 1, \texttt{Polymesher} mesh}
    \label{fig:test3:order1:polymesher}
  \end{subfigure}
  \begin{subfigure}{.49\linewidth}
    \centering
    \includegraphics[width=.8\linewidth]{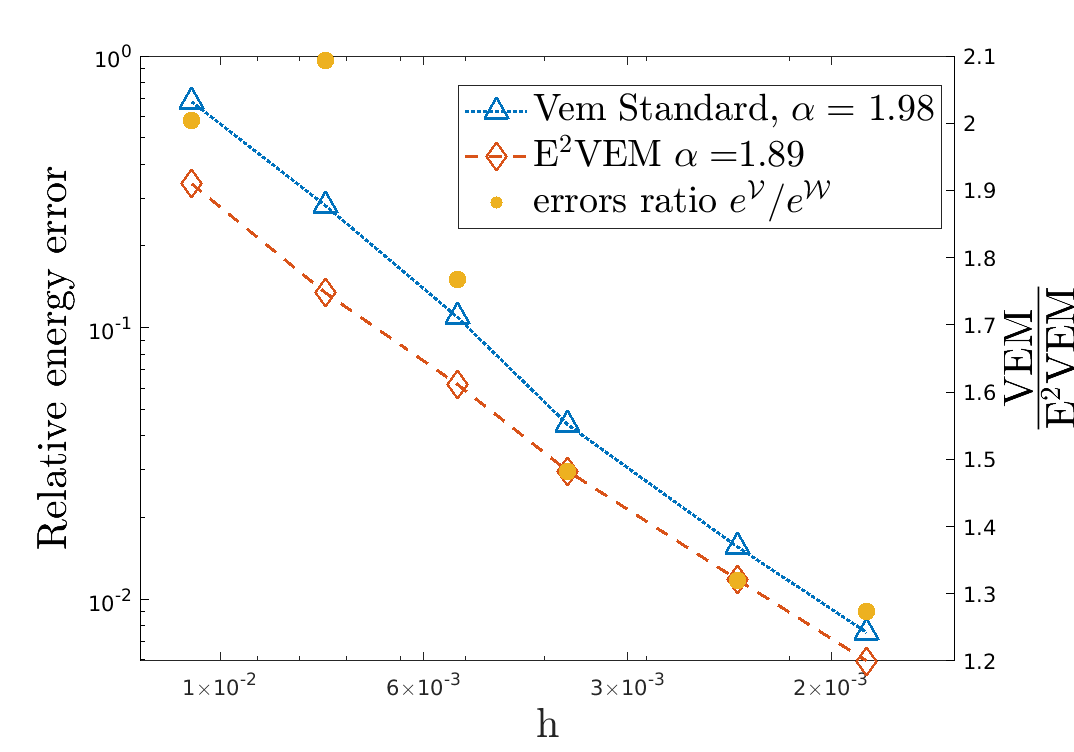}
    \caption{Order 2, \texttt{Polymesher} mesh}
  \end{subfigure}
  \caption{Test case 2. Convergence plots.}
  \label{fig:test3}
\end{figure}
\begin{table}
  \centering \input{tab/test3_projstabratio.tex}
  \caption{Test case 2. Average ratio through refinement between the infinity
    norms of the polynomial part $\mathbb{A}^\Pi$ and the stabilizing part
    $\mathbb{A}^S$ of the stiffness matrix in standard VEM.}
  \label{tab:test3_projstabratio}
\end{table}
In the second test, the exact solution is $ u(x,y) = \sin(2\pi x)\sin(80\pi y)$
and $\K = e_1 e_1^\intercal + 6.25\cdot 10^{-4} (e_2e_2^\intercal)$. In Figure
\ref{fig:test3} we display the error plots for orders 1 and 2 and Table
\ref{tab:test3_projstabratio} reports again the average ratio between the
polynomial and stabilizing parts of the standard VEM bilinear form. 
The results
are mostly consistent with the previous test, hence the convergence plots for \texttt{Cartesian} meshes are not reported for brevity. However, we observe from Figure
\ref{fig:test3:order1:polymesher} that with \texttt{Polymesher} meshes the
\EEVEM{} method reaches the asymptotic rate of convergence before the standard
VEM method. This is due to the very strong anisotropy of the solution, due to
its highly oscillating behaviour in the $y$-direction.


%% file: tab/test2_projstabratio.tex
\begin{tabular}{c|c|c||c|c}
  & \multicolumn{2}{|c||}{\texttt{Cartesian}} & \multicolumn{2}{|c}{\texttt{Polymesher}}
  \\
  \hline
  & order 1 & order 3 & order 1 & order 3
  \\
  \hline
  $\mathrm{avg}\dfrac{\norm[\infty]{\mathbb{A}^S}}{\norm[\infty]{\mathbb{A}^\Pi}}$ & 1.00 & 0.23 & 1.05 & 0.27
\end{tabular}

%% file: tab/test3_projstabratio.tex
\begin{tabular}{c|c|c||c|c}
  & \multicolumn{2}{|c||}{\texttt{Cartesian}} & \multicolumn{2}{|c}{\texttt{Polymesher}}
  \\
  \hline
  & order 1 & order 2 & order 1 & order 2
  \\
  \hline
  $\mathrm{avg}\dfrac{\norm[\infty]{\mathbb{A}^S}}{\norm[\infty]{\mathbb{A}^\Pi}}$ & 1.00 & 0.56 & 1.11 & 0.62
\end{tabular}

%% file: conclusions.tex
\section{Conclusions}
\label{sec:conclusions}

In this letter, we compared the behaviour of standard VEM and \EEVEM{} on some
Laplace test problems. Numerical results show that in the presence of strong
anisotropies of the solution and diffusivity tensor, when we apply the two
methods on general polygonal meshes, \EEVEM{} perform better than VEM in lowest
order. In all the other cases, the two methods behave equivalently.
